\documentclass[leqno,12pt]{amsart}
\usepackage{latexsym}
\usepackage{mathrsfs}
\usepackage{amsmath}
\usepackage{amssymb}
\usepackage[bookmarks]{hyperref}
\usepackage{pstricks,pst-plot}

\font\tenscr=rsfs10 
\font\sevenscr=rsfs7 
\font\fivescr=rsfs5 
\skewchar\tenscr='177 \skewchar\sevenscr='177 \skewchar\fivescr='177
\newfam\scrfam \textfont\scrfam=\tenscr \scriptfont\scrfam=\sevenscr
\scriptscriptfont\scrfam=\fivescr

\newtheorem{theorem}{Theorem}[section]
\newtheorem{lemma}[theorem]{Lemma}
\newtheorem{corollary}[theorem]{Corollary}
\newtheorem{proposition}[theorem]{Proposition}

\theoremstyle{definition}

\newtheorem{definition}[theorem]{Definition}

\allowdisplaybreaks[1]

\newcommand{\C}{\mathbb{C}}
\def\C{{\mathbf {C}\/}}

\renewcommand{\P}{\mathbb{P}}
\newcommand{\R}{\mathbf{R}}
\def\R{{\mathbf {R}\/}}

\newcommand{\range}[1]{1,\ldots, #1}

\def\CN{{\C^N}}
\def\rhull {h_r(X)}

\def\rrr#1#2{r_{#1}^{(#2)}}
\def\uuu#1#2{u_{#1}^{(#2)}}

\def\ppi#1{p_{#1}^{-1}}

\def\tk#1#2{\widetilde K_{#1}^{#2}}
\def\ck#1#2{\breve K_{#1}^{#2}}
\def\kkk#1#2{K_{#1}^{#2}}

\def\hr#1{h_r({#1})}

\def\tX{\widetilde X}
\def\wX{\widehat X}

\def\tO{\widetilde \Omega}
\def\tE{\widehat E}

\def\\{\setminus}

\newcommand{\what}{\widehat}

\newcommand{\bthm}{\begin{theorem}}
\newcommand{\ethm}{\end{theorem}}
\newcommand{\blem}{\begin{lemma}}
\newcommand{\elem}{\end{lemma}}
\newcommand{\bcor}{\begin{corollary}}
\newcommand{\ecor}{\end{corollary}}
\newcommand{\bprop}{\begin{proposition}}
\newcommand{\eprop}{\end{proposition}}
\newcommand{\bdefn}{\begin{definition}}
\newcommand{\edefn}{\end{definition}}
\newcommand{\bpf}{\begin{proof}}
\newcommand{\epf}{\end{proof}}

\newcommand{\bm}{\bibitem}

\newcommand{\bi}{\begin{itemize}}
\newcommand{\ei}{\end{itemize}}

\newcommand{\bc}{\begin{cases}}
\newcommand{\ec}{\end{cases}}

\newcommand{\ba}{\begin{array}}
\newcommand{\ea}{\end{array}}

\newcommand{\be}{\begin{equation}}
\newcommand{\ee}{\end{equation}}

\newcommand{\bea}{\begin{eqnarray}}
\newcommand{\eea}{\end{eqnarray}}

\newcommand{\beaa}{\begin{eqnarray*}}
\newcommand{\eeaa}{\end{eqnarray*}}

\newcommand{\beastar}{\begin{eqnarray*}}
\newcommand{\eeastar}{\end{eqnarray*}}

\font\tenscr=rsfs10 
\font\sevenscr=rsfs7 
\font\fivescr=rsfs5 
\skewchar\tenscr='177 \skewchar\sevenscr='177 \skewchar\fivescr='177
\newfam\scrfam \textfont\scrfam=\tenscr \scriptfont\scrfam=\sevenscr
\scriptscriptfont\scrfam=\fivescr

\def\BoX{{\mathbf B}_0(X)}
\def\Bom{{\mathbf B}_0(\Omega_m)}
\def\BX{{\mathbf B}(X)}

\def\Bone#1{{\mathbf B}(#1)}
\def\Bzero#1{{\mathbf B}_0(#1)}

\def\forallBzero#1{\hbox{for all } f\in {{\mathbf B}_0(#1)}}
\def\forallBone#1{\hbox{for all } f\in {{\mathbf B}(#1)}}

\def\textfrac#1#2{{\textstyle \frac{#1}{#2}}}

\def\ob{\overline B}
\def\od{\overline D}

\def\pfz{{{\partial f}/{\partial z_\nu}}}

\def\varep{\varepsilon}

\begin{document}
\title[Uniform algebras with dense invertible group II]{Gleason parts and point derivations\\ for uniform algebras with\\ dense invertible group II}
\author{Alexander J. Izzo}
\address{Department of Mathematics and Statistics, Bowling Green State University, Bowling Green, OH 43403}
\email{aizzo@bgsu.edu}
\thanks{The author was supported by NSF Grant DMS-1856010.}



\begin{abstract}
Due to the omission
$\vphantom{\widehat{\widehat{\widehat{\widehat{\widehat{\widehat{\wX}}}}}}}$of a hypothesis from an elementary lemma in the author's paper \lq\lq Gleason parts and point derivations for uniform algebras with dense invertible group\rq\rq, some of the proofs presented in that paper are flawed.  We prove here that nevertheless, all of the results in that paper, with the exception of the one misstated lemma, are correct.  In the process, we strengthen slightly some of the results of that paper.
\end{abstract}
\maketitle

\vskip -1.9 true in
\centerline{\footnotesize\it Dedicated to the memory of Andrew Browder} 
\vskip 1.9 truein


\section{Introduction}

In the author's paper \cite{Izzo2018}, it is shown that there exist compact sets $X$ in $\CN$ ($N\geq 2$) with nontrivial polynomial or rational hull for which the uniform algebra $P(X)$ or $R(X)$ has a dense set of invertible elements, a large Gleason part, and an abundance of bounded point derivations.  This is done by combining new constructions of \lq\lq hulls with dense invertibles\rq\rq\ with a construction of a compact set $X$ such that $R(X)$ has a Gleason part of full measure and a nonzero bounded point derivation at almost every point.  The purpose of the present paper is to correct two errors in \cite{Izzo2018} and, in particular, to show that in spite of these errors, all of the results in \cite{Izzo2018} are correct with the exception of one misstated lemma whose statement is easily corrected.

One error occurs in the proof of \cite[Theorem~1.7]{Izzo2018}.  This error is easily corrected.  See Section~\ref{Theorem17}.

The other error at first sight appears even easier to fix.  As anyone with a background in complex analysis can see, \cite[Lemma~5.1]{Izzo2018} is false as stated, but becomes true upon the addition of the hypothesis that the set $\Omega$ is connected.  The trouble is that when the lemma is invoked in \cite{Izzo2018}, the unstated hypothesis that $\Omega$ is connected is not considered.  The lemma is invoked only in the proofs of Theorems~1.5 and~6.1.  In the case of Theorem~1.5, the sets that play the role of $\Omega$ when applying the lemma are connected, as the reader can easily verify, so no change is needed in the proof.  In the case of Theorem~6.1 though, it is not at all clear that the connectedness hypothesis is satisfied.  While it is conceivable that the construction can be carried out so as to insure that the connectedness hypothesis is satisfied, we will instead give a new construction to prove \cite[Theorem~6.1]{Izzo2018}.  This new construction, which constitutes the bulk of the present paper, uses the ideas in the earlier (flawed) argument, but also involves considerable additional complications.

Note that the proofs of \cite[Theorems~1.2--1.4, and~1.6]{Izzo2018} relied upon \cite[Theorem~6.1]{Izzo2018}, so correcting the proof of \cite[Theorem~6.1]{Izzo2018} is needed to establish the validity of those results also.  (The proofs of some of these results (Theorems~1.2--1.4) do not, however, require the full force of \cite[Theorem~6.1]{Izzo2018} and could surely be established by simpler means than the argument we will give for \cite[Theorem~6.1]{Izzo2018}.)

The new construction we will give actually yields slightly more than was claimed in \cite{Izzo2018} in that the set we construct is connected.  Specifically, we will establish the following result.  Here, and throughout the paper, $B$ denotes the open unit ball $\{z: \|z\|<1\}$ in $\CN$, and $\mu$ denotes the $2N$-dimensional Lebesgue measure on $\CN$.

\bthm \label{strongcheese}
There exists a connected, compact rationally convex set $X\subset B\subset\CN$ {\rm (}$N\geq 1${\rm )} of positive $2N$-dimensional measure  such that  the collection of polynomials zero-free on $X$ is dense in $P(\ob)$  and such that there is a set $P\subset X$ of full $2N$-dimensional  measure in $X$ such that $P$ is contained in a single Gleason part for $R(X)$ and at every point of $P$ the space of bounded point derivations for 
$R(X)$ has dimension $N$.
Furthermore, given $\varep>0$, the set $X$ can be chosen so that $\mu(B \\ X)<\varep$.
\ethm

We remark that as a consequence of this theorem, the set $X$ in \cite[Theorem~1.6]{Izzo2018} can be taken to be connected.

In the next section we recall some standard definitions and notation already used above, and we introduce some other notation we will use.  Section~\ref{lemmas} is devoted to several lemmas which are used in 
Section~\ref{strongcheese-section} to prove Theorem~\ref{strongcheese} .  As mentioned earlier, the correction to the minor error in the proof of \cite[Theorem~1.7]{Izzo2018} is given in the concluding Section~\ref{Theorem17}.

It is with a mixture of joy and sorrow that I dedicate this paper to the memory of Andrew Browder.  Sorrow, of course, that he is no longer with us; joy that I had the privilege of knowing him.


\section{Preliminaries}~\label{prelim}

For
$X$ a compact Hausdorff space, we denote by $C(X)$ the algebra of all continuous complex-valued functions on $X$ with the supremum norm
$ \|f\|_{X} = \sup\{ |f(x)| : x \in X \}$.  A \emph{uniform algebra} on $X$ is a closed subalgebra of $C(X)$ that contains the constant functions and separates
the points of $X$.  

For a compact set $X$ in $\C^N$, we denote by 
$P(X)$ the uniform closure on $X$ of the polynomials in the complex coordinate functions $z_1,\ldots, z_N$, and we denote by $R(X)$ the uniform closure of the rational functions  holomorphic on (a neighborhood of) $X$.  It is well known that the maximal ideal space of $P(X)$ can be naturally identified with the \emph{polynomial hull} $\what X$ of $X$ defined by
$$\what X=\{z\in\C^N:|p(z)|\leq \max_{x\in X}|p(x)|\
\mbox{\rm{for\ all\ polynomials}}\ p
\},$$
and the maximal ideal space of $R(X)$ can be naturally identified with the \emph{rational hull} $\rhull$ of $X$ defined by
$$\rhull = \{z\in\C^N: p(z)\in p(X)\ 
\mbox{\rm{for\ all\ polynomials}}\ p
\}.$$

An equivalent formulation of the definition of $\rhull$ is that $\rhull$ consists precisely of those points $z\in \C^N$ such that every polynomial that vanishes at $z$ also has a zero on $X$.

The open unit disc in the plane will be denoted by $D$, and the open unit ball in $\CN$ will be denoted by $B$.
For $E$ a subset of $\C^N=\R^{2N}$, we will denote the $2N$-dimensional Lebesgue measure of $E$ by $\mu(E)$.
The real part of a complex number (or function) $z$ will be denoted by $\Re z$.  As in \cite{Izzo2018}, we will use the following notation in which we assume that $X$ is a compact set in $\C^N$ and that $\Omega$ is an open set in $\C^N$ each of which contains the origin:

\smallskip
\noindent
$\Bone X$ denotes the set of rational functions $f$ holomorphic on a neighborhood of $X$ such that $\|f\|_X\leq 1$,\hfil\break
$\Bzero X$ denotes the set of functions in $\Bone X$ that vanish at the origin,\hfil\break
$\Bone \Omega$ denotes the set of functions $f$ holomorphic on $\Omega$ such that \break
$\|f\|_\Omega\leq 1$, and\hfil\break
$\Bzero \Omega$ denotes the set of functions in $\Bone \Omega$ that vanish at the origin.
\smallskip

Let $A$ be a uniform algebra on a compact space $X$.
The \emph{Gleason parts} for the uniform algebra $A$ are the equivalence classes in the maximal ideal space of $A$ under the equivalence relation $\varphi\sim\psi$ if $\|\varphi-\psi\|<2$ in the norm on the dual space $A^*$.  (That this really is an equivalence relation is well-known but {\it not\/} obvious.)
For $\phi$ a multiplicative linear functional on $A$, a \emph{point derivation} on $A$ at $\phi$ is a linear functional $\psi$ on $A$ satisfying the identity 
$$\phantom{\hbox{for all\ } f,g\in A.} \psi(fg)=\psi(f)\phi(g) + \phi(f)\psi(g)\qquad \hbox{for all\ } f,g\in A.$$
A point derivation is said to be \emph{bounded} if it is bounded (continuous) as a linear functional.

We will make use of the following standard lemma.  (For a proof see \cite[Lemma~2.6.1]{Browder}.)

\begin{lemma} \label{parts}
Two multiplicative linear functionals $\phi$ and $\psi$ on a uniform algebra $A$ lie in the same Gleason part if and only if
$$\sup\{|\psi(f)|: f\in A, \|f\|\leq 1, \phi(f)=0\}< 1.$$
\end{lemma}

We will also need the following elementary observation about bounded point derivations whose easy proof is given in \cite[p.~4304]{Izzo2018}.

\begin{lemma} \label{derivations}
Let $X\subset \CN$ be a compact set.
Then the complex vector space of bounded point derivations on $R(X)$ at a particular point $x\in X$ has dimension at most $N$, and the dimension is exactly $N$ if and only if there is a number $M<\infty$ such that
\beaa
|\pfz(x)|\leq M \quad \hbox{for every } f\in \Bone {X} \hbox{\ and\ } \nu=\range N. \label{XDbound-n}
\eeaa
\end{lemma}


\section{Lemmas} \label{lemmas}
In this section we present several lemmas that will be used in the proof of Theorem~\ref{strongcheese}.

\blem\label{ratconvex}
Let $X\subset \CN$ be a compact, rationally convex set.  Then each component of $X$ is rationally convex.
\elem

It is also true that each component of a polynomially convex set is polynomially convex, and more generally, that each component of the maximal ideal space of a uniform algebra $A$ is $A$-convex, but we will not need these facts.

\bpf
Let $K$ be a component of $X$.  Since $\hr K \subset \hr X =X$, it suffices to show that given $x\in X \\ K$, there is a rational function $g$ holomorphic on a neighborhood of $X$ such that $|g(x)|> \|g\|_K$.  Choose a separation $X_0$, $X_1$ of $X$ with $K\subset X_0$ and $x\in X_1$.  By the Shilov idempotent theorem, there exists a function $f\in R(X)$ such that $f$ is identically $0$ on $X_0$ and identically $1$ on $X_1$.  Now taking for $g$ a rational function holomorphic on a neighborhood of $X$ such that $\|f-g\|_X < 1/4$ completes the proof.
\epf

\blem\label{0}
Fix $N\geq 1$ and $\varep>0$.  Let $\{p_j\}$ be a countable collection of polynomials on $\CN$ such that $p_j(0)\neq 0$ for each $j$.  Then there exists a polynomially convex Cantor set $E$ with $0\in E \subset B\subset \CN$ such that 
\item{\rm(i)} each $p_j$ is zero-free on $E$
\item{\rm(ii)} $\mu(B\\ E)<\varep$.
\elem

\bpf
The main idea of the proof is the same as the proof of \cite[Lemma~4.2]{Izzo2018}. 

Obviously it is sufficient to prove the lemma with the conditions $E\subset B$ and 
$\mu(B\\E)<\varep$ replaced by $E\subset \ob$ and $\mu(\ob \\E)<\varep$.
By multiplying each $p_j$ by a suitable complex number if necessary, we may assume that $\Re p_j(0)\neq 0$ for each $j$.  Furthermore, by enlarging the collection $\{p_j\}$, we may assume that $\{p_j\}$ is dense in $P(\ob)$.  

For each $j$, the set $p_j^{-1}(\{z\in \C: \Re z=0\})$ is a real-analytic variety in $\CN$ and hence has $2N$-dimensional measure zero.  Consequently, we can choose $0<\varep_j< \min\{ |\Re p_j(0)|, 1\}$ such that 
$$\mu\Bigl(p_j^{-1}\bigl(\{z\in \C: |\Re z|< \varep_j\}\bigl) \cap \ob \Bigr)< \varep/2^j.$$
Then
$$\mu\Bigl(\bigcup_{j=1}^\infty p_j^{-1}\bigl(\{z\in \C: |\Re z|< \varep_j\}\bigl) \cap \ob \Bigr)< \sum_{j=1}^\infty \varep/2^j = \varep.$$
Thus setting $X=\ob \\ \bigcup_{j=1}^\infty p_j^{-1}\bigl(\{z\in \C: |\Re z|< \varep_j\}\bigl)$ we have that $\mu(\ob \\ X)<\varep$.  Obviously $0\in X \subset \ob$, each $p_j$ is zero-free on $X$, and $X$ is compact.  

For each $j$, choose a closed disc $\od_j$ containing $p_j(\ob)$.  Then
$$X=\bigcap_{j=1}^\infty \Bigl(p_j^{-1}\bigl(\{z\in \od_j: |\Re z| \geq \varep_j\}\bigr) \cap \ob\Bigr).$$
Each set $\{z\in \od_j : | \Re z| \geq \varep_j\}$ is polynomially convex since it has connected complement in the plane.  Hence each set 
$p_j^{-1}\bigl(\{z\in \od_j: |\Re z| \geq \varep_j\}$ is polynomially convex (by the elementary \cite[Lemma~3.1]{Izzo2018}).  Consequently, $X$ is polynomially convex.

Let $x$ and $y$ be arbitrary distinct points of $X$.  Because $\{p_j\}$ is dense in $P(\ob)$, there is some $p_j$ such that $\Re p_j(x)>1$ and $\Re p_j(y)<-1$.  Because $| \Re p_j| \geq \varep_j$ everywhere on $X$, it follows that $x$ and $y$ lie in different components of $X$.  Consequently, $X$ is totally disconnected.

By the Cantor-Bendixson theorem \cite[Theorem~2A.1]{Mos} the compact metric space $X$ is the disjoint union of a perfect set $E$ and an at most countable set.  Then $\mu(\ob \\ E)=\mu(\ob \\ X)$, and $E$ is a Cantor set by the usual characterization of Cantor sets as the compact, totally disconnected, metrizable spaces without isolated points.  Finally, $E$ is polynomially convex, for $\tE$ is contained in $X$, so if $\tE$ were strictly larger than $E$, then $\tE \\ E$ would have isolated points, which is impossible (by the Shilov idempotent theorem, for instance).
\epf

\blem\label{A}
Fix $N\geq 1$ and $\varep>0$.  Then there exists a compact rationally convex set $Y$ with $0\in Y \subset B\subset \CN$ such that 
\item{\rm(i)} the set of polynomials zero-free on $Y$ is dense in $P(\ob)$
\item{\rm(ii)} $\mu(B\\ Y)<\varep$
\item{\rm(iii)} $Y$ is connected.
\elem

\bpf
Choose a countable collection $\{p_j\}$ of polynomials that is dense in $P(\ob)$ and such that $p_j(0)\neq 0$ for each $j$.  By the preceding lemma, there exists a polynomially convex Cantor set $E$ such that $0\in E \in B$, each $p_j$ is zero-free on $E$, and $\mu(B\\ E) < \varep$.  By Antoine's theorem \cite{Antoine} (or see \cite{Whyburn} where a more general result is given), there is an arc $C$ in $\CN$ that contains $E$.  Let $\sigma: [0,1]\rightarrow C$ be a homeomorphism.  Assume without loss of generality that $\sigma(0)$ and $\sigma(1)$ are in the Cantor set $E$, and let $K$ be the Cantor set $\sigma^{-1}(E)$.  Let $(a_1,b_1)$, $(a_2,b_2)$, $\ldots$ be the disjoint open intervals whose union is $[0,1]\\ K$.  Define a map $\gamma: [0,1] \rightarrow \CN$ by setting $\gamma=\sigma$ on $K$ and taking $\gamma$ to be affine on each interval $[a_k,b_k]$ for $k=1, 2, \ldots$.  The reader can verify that $\gamma$ is continuous.  Consequently, $\gamma([0,1])$ is a connected, compact subset of $B \subset \CN$ that contains $E$.  Set $J=\gamma([0,1]$.  Set $I_k=[a_k, b_k]$ for $k=1,2,\ldots$.  Because each $p_j$ is zero-free on $E$, there exists $\varep_j>0$ such that $p_j(E)$ is disjoint from $\{z\in \C: |z|<\varep_j\}$.  Each set $p_j(I_k)$ has empty interior in the plane, so it follows that the set $p_j(J)\cap \{z\in \C: |z|<\varep_j\}$ is a countable union of nowhere dense set and thus has empty interior.  Thus there exist arbitrarily small complex numbers $\alpha$ such that $p_j+\alpha$ has no zeros on $J$.  Consequently, conditions~(i), (ii), and~(iii) all hold with $Y$ replaced by $J$.  

Now set $Y=\hr J$.  Then conditions (i) and (ii) are immediate, and condition (iii) follows from the connectedness of $J$ by Lemma~\ref{ratconvex}.
\epf

The proof of the next lemma is an easy exercise.   
However, given that the incorrect statement of the lemma in \cite{Izzo2018} is what created the need for the present paper, and given the central role the lemma places in the proof of our main result, we include the proof in full.

\blem\label{5.1}
Suppose $K$ is a compact set containing the origin and contained in an open set $\Omega\subset \CN$.  Then there exists an $R$ with $0<R<1$ such that $\|f\|_K\leq R$ for all $f\in \Bzero \Omega$ if and only if $K$ is contained in the component of $\Omega$ that contains the origin.
\elem

\bpf
Suppose that $K$ is contained in the component of $\Omega$ that contains the origin.  Then we may assume without loss of generality that $\Omega$ is connected.  Now assume to get a contradiction that no $R$ as in the statement of the lemma exists.  Then for each $n=1, 2,\ldots$, there exists a function $f_n\in \Bzero \Omega$ such that $\|f_n\|_K \geq 1 - 1/n$.  By Montel's theorem, $\Bzero \Omega$ is a normal family, so some subsequence of $(f_n)$ converges uniformly on compact subsets of 
$\Omega$ to a limit function $f$.  Then $f$ is a holomorphic function on $\Omega$ such that $f(0)=0$, $\|f\|_\Omega \leq 1$, and $\|f\|_K \geq 1$.  But then $f$ is a nonconstant holomorphic function on the connected, open set $\Omega$ that takes on a maximum in modulus at some point, a contradiction.

Now suppose that $K$ fails to be contained in the component of $\Omega$ that contains the origin.  Let $g$ be the function on $\Omega$ that is identically $0$ on the component of $\Omega$ that contains the origin and is identically $1$ everywhere else on $\Omega$.  Then $\|g\|_K=1$ and $g\in \Bzero \Omega$.  Thus no $R$ as in the statement of the lemma exists.
\epf

The next three lemmas are taken from \cite{Izzo2018}.

\blem\cite[Lemma~5.2]{Izzo2018}\label{5.2}
Suppose $K$ is a compact set contained in an open set $\Omega\subset \CN$.  Then there exists an $M<\infty$  such that $\|\pfz\|_K\leq M$ for all $f\in \Bone \Omega$ and $\nu=\range N$.
\elem

\blem\cite[Lemma~6.2]{Izzo2018}\label{6.2}
Let $K$ be a compact set containing the origin and contained in an open set $\Omega\subset \CN$, let $p$ be a polynomial on $\CN$ with no zeros on $K$, and let $\varep>0$ be given.  Let $0<R<1$, and suppose that $\|f\|_K\leq R$ for all $f\in \Bzero \Omega$.
Then there exists an $r>0$ such that $p^{-1}(r\od)\cap K=\emptyset$ and $\|f\|_K\leq R+\varep$ for all $f\in \Bzero {\Omega\\ p^{-1}(r\od)}$.
\elem

\blem\cite[Lemma~6.3]{Izzo2018}\label{6.3}
Let $K$ be a compact set contained in an open set $\Omega\subset \CN$, let $p$ be a polynomial on $\CN$ with no zeros on $K$, and let $\varep>0$ be given.  Let $0<M<\infty$, and suppose that $\| \pfz \|_K\leq M$ for all $f\in \Bone \Omega$ and $\nu=\range N$.
Then there exists an $r>0$ such that $p^{-1}(r\od) \cap K=\emptyset$ and $\|\pfz\|_K\leq M+\varep$ for all $f\in \Bone {\Omega\\ p^{-1}(r\od)}$ and $\nu=\range N$.
\elem

\blem\label{shrink1}
Let $K$ be a compact set containing the origin and contained in an open set 
$\Omega\subset \CN$, and let $\varep>0$ be given.  Let $0<R<1$, and suppose that $\|f\|_K\leq R$ for all $f\in \Bzero \Omega$.
Let $(\Omega_n)$ be an increasing sequence of open sets of $\CN$ with union $\Omega$.
Then there exists an $n$  such that $\Omega_n\supset K$ and $\|f\|_K\leq R+\varep$ for all $f\in \Bzero {\Omega_n}$.
\elem

\bpf
Assume to get a contradiction that no such $n$ exists.  Let $J$ be a positive integer large enough that $\Omega_n \supset K$.  Then for each $n=J, J+1, \ldots$, there exists a function $f_n\in \Bzero {\Omega_n}$ such that $\|f_n\|_K> R+\varep$.  For each $m=J, J+1, \ldots$, the set $\{f_n: n\geq m\}$ is a normal family on $\Omega_m$.
Thus there is a subsequence of $(f_n)$ that converges uniformly on compact subsets of $\Omega_J$.  We may then choose a further subsequence that converges uniformly on compact subsets of  $\Omega_{J+1}$.  Continuing in this manner taking subsequences of subsequences, and then applying the usual diagonalization argument, we arrive at a subsequence $(f_{n_k})$ of $(f_n)$ such that for each compact subset $L$ of $\Omega$, the sequence  $(f_{n_k})$
converges uniformly on $L$.  Note that given $L$, there may be a finite number of terms of the sequence $(f_{n_k})$ that are not defined on $L$, but this does not matter.  Thus there is a well-defined limit function $f$ that is holomorphic on 
$\Omega$.  The function $f$ is in $\Bzero\Omega$.  Thus by hypothesis, $\|f\|_K\leq R$.  But  since $f_{n_k}\rightarrow f$ uniformly on $K$, and $\|f_n\|_K>R+\varep$ for all $n$, this is a contradiction.
\epf

\blem\label{shrink2}
Let $K$ be a compact set contained in an open set $\Omega\subset \CN$, and let $\varep>0$ be given.  Let $0<M<\infty$, and suppose that $\| \pfz \|_K\leq M$ for all $f\in \Bone \Omega$ and $\nu=\range N$.
Let $(\Omega_n)$ be an increasing sequence of open sets of $\CN$ with union $\Omega$.
Then there exists an $n$ such that $\Omega_n\supset K$ and $\|\pfz\|_K\leq M+\varep$ for all $f\in \Bone {\Omega_n}$ and $\nu=\range N$.
\elem

\bpf
This is proven by a normal families argument similar to the one just presented and hence is left to the reader.
\epf

\blem\label{shrink0}
Given $K\subset U \subset \CN$ with $K$ a compact set and $U$ a connected open set, there exists a connected open set $W$ with compact closure such that $
K\subset W \subset \overline W \subset U$.
\elem

\bpf
Cover $K$ by a finite number of open balls $B_1, \ldots, B_k$ whose closures 
are contained in $U$.  Let $a_1, \ldots, a_k$ be the centers of these balls.  For each $j=2,\ldots, k$, choose a path $\gamma_j$ from $a_1$ to $a_j$ in $U$.  Let $\gamma_j^*$ denote the image of the path $\gamma_j$.  Set $J=\ob_1 \cup \cdots \cup \ob_k \cup \gamma_2^* \cup \cdots \cup \gamma_k^*$.  Then $J$ is a compact, connected subset of $U$ that contains $K$.  Cover $J$ by a finite number of open balls $W_1, \ldots, W_m$ whose closures are contained in $U$.  Assume without loss of generality that each $W_j$ intersects $J$.  Set $W=W_1\cup \cdots \cup W_m$.  Then $W$ is as asserted in the statement of the lemma.
\epf


\section{Proof of Theorem~\ref{strongcheese}} \label{strongcheese-section}

Given a real number $r>0$ and a subset $S$ of $\CN$, we will denote by $rS$ the set $\{ rs: s\in S\}$.
Recall the notations
$\Bone X$,
$\Bzero X$,
$\Bone \Omega$, and
$\Bzero \Omega$ introduced in Section~\ref{prelim}.

\bpf[Proof of Theorem~\ref{strongcheese}]
Fix $\varep>0$.  Let $Y$ be the set whose existence is given by Lemma~\ref{A}, and  choose a countable collection $\{p_j\}$ of polynomials on $\CN$ that is dense in $P(\ob)$ and such that each $p_j$ is zero-free on $Y$. 

It is to be understood that throughout the proof, all subscripts and superscripts are positive integers.
We will choose a strictly decreasing sequence $(s_m)_{m=1}^\infty$ of strictly positive numbers less than 1
and a doubly indexed collection $\{\rrr jm:1\leq j\leq m\}$ of strictly positive numbers such that for each $j$, the sequence $(\rrr jm)_{m=j}^\infty$ is strictly increasing and bounded.  Given these, we define $\tX_m$ and $\tO_m$ by
\be
\tX_m=s_m\ob\\ \bigcup_{j=1}^m p_j^{-1}(\rrr jm D) \label{txm}
\ee
and
\be
\tO_m=s_m B\\ \bigcup_{j=1}^m p_j^{-1}(\rrr jm \od). \label{tom}
\ee
The numbers $\rrr jm$ will be shown small enough that $\Omega_m$, and hence also $X_m$, contains the origin. 
We denote the component of each of $\tX_m$ and $\tO_m$ that contains the origin respectively by $X_m$ and $\Omega_m$.  Note that $X_m$ is compact and $\Omega_m$ is open.  Note also that 
$$\tX_1\supset \tO_1\supset \tX_2\supset \tO_2\supset \tX_3\supset \tO_3\supset \cdots$$
and
$$X_1\supset \Omega_1\supset X_2\supset \Omega_2\supset X_3\supset \Omega_3\supset \cdots.$$

Let 
$$X=\bigcap_{m=1}^\infty X_m.$$

Letting $s=\inf s_m$ and $r_j=\sup\limits_m \rrr jm$, we have that the intersection $\tX=\bigcap_{m=1}^\infty \tX_m$ is given by
$$\tX= s\ob \\ \bigcup_{j=1}^\infty p_j^{-1}(r_j D),$$
and since the intersection of a decreasing sequence of compact connected sets is connected, $X$ is the component of $\tX$ that contains the origin.

The choice of $(s_m)$ and $\{\rrr jm\}$ will involve also choosing an increasing sequence $(t_n)$ of strictly positive numbers and a doubly indexed collection $\{\uuu jn: j=1, 2,\ldots\ {\rm and\ } n=1,2,\ldots\}$ of strictly positive numbers such that for each $j$, the sequence $(\uuu jn)_{n=1}^\infty$ is decreasing.  The choices will be made in such a way that $t_\alpha < s_\beta$ and $\uuu j\beta  > \rrr j\alpha$ for all $\alpha$, $\beta$, and $j$ for which the quantities are defined.
We then define $\tk nm$ for all $m\geq n\geq 1$ by
\be
\tk nm=t_n \ob\\ \bigcup_{j=1}^m p_j^{-1}(\uuu jn D). \label{K}
\ee
Note that then $\tk nm$ is contained in $\tO_m$.
We set
\be
\kkk nm=\tk nm\cap \Omega_m. \label{knm}
\ee
As a component of $\tO_m$, the set $\Omega_m$ is closed in $\tO_m$, and consequently, $\kkk nm$ is compact.

We will show that $(s_m)$, $\{\rrr jm\}$, $(t_n)$, and $\{\uuu jn\}$ can be chosen so that the following conditions hold for all $1\leq n\leq m$ and $1\leq j\leq m$.  (We present our required  conditions in a slightly odd way so as to set up a proof by induction.  The first seven of these conditions have already been mentioned above.)

{\baselineskip=18pt
\begin{enumerate}
\item[(i)] The numbers $s_m$, $\rrr jm$, $t_n$, and $\uuu jn$ are strictly positive.
\item[(ii)] $1>s_1 > s_2 > \cdots > s_m$
\item[(iii)]  $r_j^{(j)} < r_j^{(j+1)} < \cdots < \rrr jm$
\item[(iv)] $t_1\leq t_2\leq \cdots \leq t_m$
\item[(v)]  $\uuu j1 \geq \uuu j2 \geq \cdots \geq \uuu jm$
\item[(vi)]  $t_m<s_m$ (and hence, $t_\alpha< s_\beta$ for all $\alpha$ and $\beta$)
\item[(vii)]  $\uuu jm > \rrr jm$ (and hence $\uuu j\beta > \rrr j\alpha$  for all $j \leq \alpha \leq m$ and $1 \leq \beta \leq m$)
\item[(viii)]  $s_mB\supset Y$
\item[(ix)]  $\ppi j (\rrr jm \od)$ is disjoint from $Y$
\item[(x)]  $0\notin \ppi j (\uuu jn D)$
\item[(xi)] $\mu(X_n\\ K_n^n) \leq (1/n) \mu(B)$
\item[(xii)]  $\mu(\ppi j (\uuu jn D) \cap B) \leq (1/2^{j-1}) \mu(B)$
\item[(xiii)]  There exists $0<R_n<1$ such that 
$$\|f\|_{K_n^m}\leq R_n +(\textfrac 14 + \textfrac 18 + \cdots + \textfrac 1{2^{m-n+1}})(1-R_n) \quad \forallBzero {\Omega_m}.$$
(When $m=n$, the sum in parentheses is undertood to be zero.)
\item[(xiv)]   There exists $0<M_n<\infty$ such that 
$$\|\pfz\|_{K_n^m}\leq (1+ \textfrac 12 +\cdots + \textfrac 1{2^{m-n}})M_n \quad \forallBone {\Omega_m} \hbox{\ and\ } \nu=\range N.$$
\end{enumerate}
}

Assuming for the moment that $(s_m)$, $\{\rrr jm\}$, $(t_n)$, and $\{\uuu jn\}$ have been chosen so that the above conditions are satisfied, we now prove that the set $X$ has the properties stated in the theorem.  Clearly $\tX$ is rationally convex, and hence the component $X$ of $\tX$ is rationally convex by Lemma~\ref{ratconvex}.  By conditions~(viii) and~(ix), $\tX\supset Y$.  Therefore, $X\supset Y$ by the connectedness of $Y$, and hence $\mu(B \\ X)<\varep$.  Since each $p_j$ is zero-free on $X$, the collection of polynomials zero-free on $X$ is dense in $P(\ob)$.  Finally, for each $n=1,2,\ldots$, set $K_n=\bigcap_{m=n}^\infty \kkk nm$, and set $P=\bigcup_{n=1}^\infty K_n$.  Then $P$ is contained in $X$, and we now show that $P$ is as asserted in the statement of the theorem.

Note that for every $n$,
$$X\\ P \subset X_n\\ K_n \subset
\Big[ X_n \\ K_n^n\Big] \bigcup \left[\, \bigcup_{j=n+1}^\infty (\ppi j (\uuu jn D) \cap \ob)\right].$$
Letting $n\rightarrow \infty$, we see that conditions (xi) and (xii) imply that $\mu(X\\P)=0$.

Now consider an arbitrary point $z\in P$.  Fix $n$ such that $z\in K_n$.  Given $g\in \BoX$, there is some $m\geq n$ such that $g\in \Bom$.  Since $z\in \kkk nm$, condition~(xiii) gives that 
$$|g(z)|\leq R_n +(\textfrac 14 + \textfrac 18 + \cdots + \textfrac 1{2^{m-n+1}})(1-R_n).$$
Therefore,
$$\sup_{f\in \BoX} |f(z)| \leq R_n + \textfrac 12 (1-R_n) <1.$$
Thus by Lemma~\ref{parts}, $z$ lies in the Gleason part of the origin for $R(X)$.  Similarly, using condition~(xiv), we get that 
$$\sup_{f\in \BX} |(\pfz)(z)| \leq 2M_n<\infty \hbox{\ for\ all\ } \nu=\range N.$$
Thus by Lemma~\ref{derivations}, the space of bounded point derivations on $R(X)$ at $z$ has dimension $N$.

It remains to choose $(s_m)$, $\{\rrr jm\}$, $(t_n)$, and $\{\uuu jn\}$.  First choose $0<s_1< 1$ large enough that $s_1 B\supset Y$, and choose $r_1^{(1)} >0$ small enough that $p_1^{-1}(r_1^{(1)} \od)$ is disjoint from $Y$.  Then choose $0<t_1<s_1$ and $u_1^{(1)} > r_1^{(1)}$ with $u_1^{(1)}$ close enough to $r_1^{(1)}$ that $0\notin p_1^{-1}(u_1^{(1)} D)$.  
The sets $\tX_1$, $\tO_1$, $X_1$, $\Omega_1$, $\tk 11$, and $\kkk 11$ are now defined in accordance with equations~(\ref{txm}), (\ref{tom}), (\ref{K}) and~(\ref{knm}), and the remarks that immediately follow them.
By Lemmas~\ref{5.1} and~\ref{5.2}, there exist $0< R_1 <1$ and $0< M_1 <\infty$ such that
$$\|f\|_{K_1^1}\leq R_1 \quad \forallBzero {\Omega_1}$$
and
$$\|\pfz\|_{K_1^1}\leq M_1 \quad \forallBone {\Omega_1} \hbox{\ and\ } \nu=\range N.$$
Thus conditions~(i)--(xiv) hold for $n=m=j=1$.

We now continue by induction.  Suppose that for some $k\geq 1$ we have chosen $s_m$, $\rrr jm$, $t_n$, and $\uuu jn$ for all $1\leq n\leq m\leq k$ and $1\leq j \leq m \leq k$ such that conditions~(i)--(xiv) hold for all these values of $n$, $m$, and $j$.  We will refer to conditions~(i)--(xiv) for these values of $n$, $m$, and $j$ as conditions~(i)--(xiv) of the induction hypothesis.

Choose $\uuu {k+1}1=\uuu{k+1}2=\cdots=\uuu{k+1}k>0$ small enough that for each $n=1,\ldots, k$ we have that $0\notin \ppi {k+1} (\uuu {k+1}n D)$ and $\mu(\ppi {k+1} (\uuu {k+1}n D) \cap B)< (1/2^k) \mu(B)$.  Equation (\ref{K}) now defines sets $\tk 1{k+1},\ldots, \tk k{k+1}$.
For $n=1,\ldots, k$, set 
$$\ck n{k+1}= \tk n{k+1} \cap \Omega_k,$$
and note that then $\ck n{k+1} \subset \kkk nk$.  Note also that each of $\ck 1{k+1},\ldots, \ck k{k+1}$ is compact.  Lemmas~\ref{6.2} and~\ref{6.3}, together with conditions~(xiii) and~(xiv) of the induction hypothesis, yield the existence of $0<c<\uuu {k+1}k$ such that for all $n=1,\ldots, k$, we have
\begin{equation}\openup1\jot
\begin{split}
\|f\|_{\ck n{k+1}}&\leq R_n +(\textfrac 14 + \textfrac 18 + \cdots + \textfrac 1{2^{k-n+1}} + \textfrac {1/2}{2^{(k+1)-n+1}})(1-R_n)\cr
& \hskip 120pt \forallBzero {\Omega_k \setminus \ppi {k+1} (c\od)}
\end{split} \label{prebound1}
\end{equation}
and
\begin{equation}\openup1\jot
\begin{split}
\|\pfz\|_{\ck n{k+1}}&\leq (1+ \textfrac 12 +\cdots + \textfrac 1{2^{k-n}} + \textfrac {1/2}{2^{(k+1)-n}})M_n\cr 
&\forallBone {\Omega_k \setminus \ppi {k+1} (c\od)} \hbox{\ and\ } \nu=\range N.
\end{split} \label{prebound2}
\end{equation}
Note that by choosing $c$ small enough, we can arrange to have $\ppi {k+1}(c\od)$ disjoint from $Y$.

Denote the component of $\Omega_k\\ \ppi {k+1} (c\od)$ containing the origin by $U$.
By the regularity of Lebesgue measure, there exists a compact set $K$ contained in $U$ such that $\mu(U\\ K)< \textfrac {1}{2} (1/(k+1))\mu(B)$.  By Lemma~\ref{shrink0}, there exists a connected open set $W$ such that $K\subset W\subset \overline W \subset U$.  Then of course
\be
\mu(U\\ \overline W) <  \textfrac {1}{2} (1/(k+1))\mu(B). \label{measure}
\ee
If we let $(\alpha_n)$, $(\beta_j^{(n)})_{n=1}^\infty$ for $j=1,\ldots, k$, and $(\gamma_n)$ be sequences of strictly positive numbers with $(\alpha_n)$ increasing to $s_k$, with $(\beta_j^{(n)})_{n=1}^\infty$ decreasing to $r_j^{(k)}$, and with $(\gamma_n)$ decreasing to $c$, and we set 
$$V_n=\alpha_n B \\ \left [\left(\, \bigcup_{j=1}^k \ppi j (\beta_j^{(n)} \od) \right ) \cup \ppi {k+1} (\gamma_n \od) \right ],$$
then $(V_n)$ is an increasing sequence of open sets with union $\tO_k \\ \ppi {k+1} (c\od)$.  Therefore, by Lemmas~\ref{shrink1} and~\ref{shrink2}, it follows from inequalities~(\ref{prebound1}) and~(\ref{prebound2}) that there exist numbers $s_{k+1}$ and $\rrr j{k+1}$ for $j=1,\ldots, k+1$ such that $0<t_k< s_{k+1}< s_k$, such that $\rrr jk < \rrr j{k+1} < \uuu jk$ for $j=1,\ldots, k$, such that $c < \rrr{k+1}{k+1} < \uuu {k+1}k$, and such that 
defining $\tO_{k+1}$ in accordance with equation~(\ref{tom}), we have that $\tO_{k+1}$ contains each of $\ck 1{k+1},\ldots, \ck k{k+1}$ and 
we have, for each $n=1,\ldots k$, the inequalities 
\begin{equation}\openup1\jot  \label{b1}
\begin{split}
\|f\|_{\ck n{k+1}}&\leq R_n +(\textfrac 14 + \textfrac 18 + \cdots + \textfrac 1{2^{k-n+1}} + \textfrac 1{2^{(k+1)-n+1}})(1-R_n)\cr
& \hskip 165pt \forallBzero {\tO_{k+1}} 
\end{split}
\end{equation}
and
\begin{equation}\openup1\jot  \label{b2}
\begin{split} 
\|\pfz\|_{\ck n{k+1}}&\leq (1+ \textfrac 12 +\cdots + \textfrac 1{2^{k-n}} + \textfrac 1{2^{(k+1)-n}})M_n\cr 
&  \hskip 45pt \forallBone {\tO_{k+1}} \hbox{\ and\ } \nu=\range N.
\end{split}
\end{equation}
Then by Lemma~\ref{5.1}, each of  $\ck 1{k+1},\ldots, \ck k{k+1}$ is contained in the component $\Omega_{k+1}$ of $\tO_{k+1}$ containing the origin, and  inequalities~(\ref{b1}) and~(\ref{b2}) continue to hold with $\tO_{k+1}$ replaced by $\Omega_{k+1}$.  Furthermore, by the the compactness of $Y$ and $\overline W$, the numbers $s_{k+1}$ and $\rrr 1{k+1},\ldots, \rrr {k+1}{k+1}$ can be chosen so that in addition, $s_{k+1}B\supset Y$ and for $j=1,\ldots, k+1$, the set $\ppi j (\rrr j{k+1} \od)$ is disjoint from $Y$ (and hence $Y \subset \tO_{k+1} $), and also 
$\overline W\subset \tO_{k+1}$.  Then by the connectedness of $Y$ and $\overline W$, we have that $Y\subset \Omega_{k+1}$ and $\overline W\subset \Omega_{k+1}$.  

With $s_{k+1}$ and $\rrr 1{k+1},\ldots, \rrr {k+1}{k+1}$  now chosen, $\tX_{k+1}$ is defined in accordance with equation~(\ref{txm}), and $X_{k+1}$ is defined to be the component of $\tX_{k+1}$ that contains the origin.  
Note that $\tX_{k+1} \subset \Omega_k \\ \ppi {k+1} (c\od)$, and so 
$X_{k+1} \subset U$.  Consequently, $X_{k+1}\\ \Omega_{k+1} \subset U\\ \overline W$.  Therefore, by inequality (\ref{measure}) we have
\be\label{mostmeasure}
\mu(X_{k+1} \\ \Omega_{k+1}) \leq \mu (U\\ \overline W) < \textfrac {1}{2} (1/(k+1)) \mu(B).
\ee
We define the sets $\kkk 1{k+1},\ldots, \kkk k{k+1}$ in accordance with equation~(\ref{knm}).  Observe that then
$$\kkk n{k+1} = \ck n{k+1} \quad \hbox{for all } j=1,\ldots, k.$$
Therefore, since inequalities~(\ref{b1}) and~(\ref{b2}) hold with $\tO_{k+1}$ replaced by $\Omega_{k+1}$, we have for all $n=1,\ldots, k$ that
\begin{equation}\nonumber\openup1\jot
\begin{split}
\|f\|_{\kkk n{k+1}}&\leq R_n +(\textfrac 14 + \textfrac 18 + \cdots + \textfrac 1{2^{k-n+1}} +\textfrac 1{2^{(k+1)-n+1}})(1-R_n)\cr
& \hskip 165pt\forallBzero {\Omega_{k+1}}
\end{split} 
\end{equation}
and
\begin{equation}\nonumber\openup1\jot
\begin{split}
\|\pfz\|_{\kkk n{k+1}}&\leq (1+ \textfrac 12 +\cdots + \textfrac 1{2^{k-n}} + \textfrac 1{2^{(k+1)-n}})M_n \cr
& \hskip 45pt \forallBone {\Omega_{k+1}} \hbox{\ and\ } \nu=\range N.
\end{split} 
\end{equation}

Finally, choose $t_{k+1}$ and $\uuu 1{k+1}, \ldots, \uuu {k+1}{k+1}$ such that $0<t_k\leq t_{k+1} < s_{k+1}$ and 
$\uuu jk \geq \uuu j{k+1} >\rrr j{k+1}$ for $j=1,\ldots k+1$ with $t_{k+1}$ close enough to $s_{k+1}$ and each $\uuu j{k+1}$ close enough to $\rrr j{k+1}$ that
$0\notin \ppi j(\uuu j{k+1}D)$ and defining $\tk {k+1}{k+1}$ in accordance with equation~(\ref{K}) we have that $\tk {k+1}{k+1} \subset \tO_{k+1}$ and 
$$\mu(\tO_{k+1} \\ \tk {k+1}{k+1}) < \textfrac {1}{2} (1/(k+1)) \mu(B).$$
Then with $\kkk {k+1}{k+1}$ defined in accordance with equation~(\ref{knm}),
we also have 
$$\mu(\Omega_{k+1} \\ \kkk {k+1}{k+1}) < \textfrac {1}{2} (1/(k+1)) \mu(B).$$
Combined with inequation~(\ref{mostmeasure}), this gives that
$$\mu(X_{k+1} \\ \kkk {k+1}{k+1}) < (1/(k+1)) \mu(B).$$

Lemmas~\ref{5.1} and~\ref{5.2} yield numbers $0< R_{k+1} <1$ and $0< M_{k+1} < \infty$ such that
$$
\|f\|_{\kkk {k+1}{k+1}}\leq R_{k+1}  \quad \forallBzero {\Omega_{k+1}}
$$
and
$$
\|\pfz\|_{\kkk {k+1}{k+1}}\leq M_{k+1} \quad \forallBone {\Omega_{k+1}} \hbox{\ and\ } \nu=\range N.
$$
Now, as the reader can verify, we have chosen $s_m$, $\rrr jm$, $t_n$, and $\uuu jn$ for all $1\leq n \leq m\leq k+1$ and $1\leq j\leq m \leq k+1$ in such a way that conditions (i)--(xiv) hold for all these values of $n$, $m$, and $j$.  Thus the induction can proceed, and the proof of Theorem~\ref{strongcheese} is complete.
\epf


\section{Correction to the proof of \cite[Theorem~1.7]{Izzo2018}}~\label{Theorem17}

In the proof of \cite[Theorem~1.7]{Izzo2018} there are compact sets $K\subset X\subset \CN$, an injective holomorphic map 
$\sigma:D \rightarrow \wX \\ K$, and a polynomial $p$ on $\CN$ such that $p(\sigma(0))=0$ and $p$ has no zeros on $K$.  It is asserted that \lq\lq because 
$\sigma$ is injective, we may assume, by adding to $p$ a small multiple of a suitable first degree polynomial if necessary,
that  $\bigl(\partial (p\circ \sigma)/\partial z\bigr)(0)\not= 0$".  This is of course false since it could be that $\bigl(\partial \sigma/\partial z\bigr)(0)= 0$.  (For instance, $\sigma$ could be given by $\sigma(z)=(z^2, z^3)$.)  However, all that is needed in the proof of \cite[Theorem~1.7]{Izzo2018} is the ability to assume that $p$ is not identically zero on $\sigma(D)$.  This can be achieved 
by adding to $p$, if necessary, a small multiple of a polynomial that vanishes at $\sigma(0)$ but is not identically zero on 
$\sigma(D)$.

A generalization of \cite[Theorem~1.7]{Izzo2018}, whose proof runs along the same lines, is given in the author's paper \cite{Izzo-spaces}.  There the argument is presented correctly; see the proof of \cite[Theorem~5.1]{Izzo-spaces}.

\end{document}